\documentclass[11pt]{amsart}

\usepackage{amssymb,amsmath,amsthm,newlfont}
\usepackage[shortlabels]{enumitem}
\usepackage{hyperref}

\theoremstyle{plain}
\newtheorem{theorem}{Theorem}[section]
\newtheorem{proposition}[theorem]{Proposition}

\newtheorem{lemma}[theorem]{Lemma}

\theoremstyle{definition}
\newtheorem{definition}[theorem]{Definition}

\theoremstyle{remark}

\numberwithin{equation}{section}

\newcommand{\DD}{\mathbb{D}}

\newcommand{\cH}{\mathcal{H}}

\newcommand{\mC}{\mathbb{C}}
\newcommand{\mR}{\mathbb{R}}
\newcommand{\mN}{\mathbb{N}}

\newcommand{\mD}{\mathbb{D}}

\newcommand{\dimH}{\dim_{\mathrm H}}
\newcommand{\dimP}{\dim_{\mathrm P}}
\newcommand{\dimMl}{\underline{\dim}_{\mathrm M}}
\newcommand{\dimMu}{\overline{\dim}_{\mathrm M}}

\begin{document}

\date{14 October 2026}

\title[Holomorphic motion]{A  holomorphic motion with non-subharmonic Hausdorff dimension}

\author[T. Ransford]{Thomas Ransford}
\address{D\'epartement de math\'ematiques et de statistique, Universit\'e Laval,
Qu\'ebec City (Qu\'ebec),  Canada G1V 0A6.}
\email{ransford@mat.ulaval.ca}

\author[M. Younsi]{Malik Younsi}
\address{Department of Mathematics, University of Hawaii Manoa, Honolulu, HI 96822, USA.}
\email{malik.younsi@gmail.com}

\thanks{Ransford supported by NSERC Discovery Grant RGPIN-2026-04565. Younsi supported by NSF Grant DMS-2050113.}

\begin{abstract}
We present an example of a holomorphic motion of a Cantor set  $\lambda \mapsto A_\lambda$ 
for which neither the Hausdorff dimension $\dimH(A_\lambda)$ nor the lower Minkowski dimension
$\dimMl(A_\lambda)$ is a subharmonic function of~$\lambda$.
This answers negatively a question of Fuhrer--Ransford--Younsi,
and distinguishes these dimensions from packing and upper Minkowski dimensions, 
for which subharmonicity holds. All four dimensions are computed explicitly.
\end{abstract}


\makeatletter
\@namedef{subjclassname@2020}{\textup{2020} Mathematics Subject Classification}
\makeatother


\maketitle

\section{Introduction}\label{S:intro}

In what follows, we write $\mD$ and $\mC$ for the open unit disk and the complex plane, respectively.

\begin{definition}\label{D:holomotion}
Let $A$ be a subset of  $\mC$. A \textit{holomorphic motion} of $A$ is a map $f:\mD\times A \to \mC$ such that:
\begin{enumerate}[(i)]
\item for each fixed $z\in A$, the map $\lambda\mapsto f(\lambda,z)$ is holomorphic on $\mD$;
\item for each fixed $\lambda \in\mD$, the map $z\mapsto f(\lambda,z)$ is injective on $A$;
\item for all $z\in A$, we have $f(0,z)=z$.
\end{enumerate}
We write $f_\lambda(z):=f(\lambda,z)$ and $A_\lambda:=f_\lambda(A)$.
We sometimes abuse terminology and call the set-valued map $\lambda\mapsto A_\lambda$ a holomorphic motion of~$A$.
\end{definition}

\begin{definition}\label{D:infharmonic}
A positive function $u:\mD\to[0,\infty)$ is called
\emph{inf-harmonic} if there exists a  family $\cH$ of harmonic functions on $\mD$ such that $u(\lambda)=\inf_{h\in\cH}h(\lambda)$
for all $\lambda\in \mD$.
\end{definition}
For more details on holomorphic motions and inf-harmonic functions, we refer the reader to \cite{FRY}. In that paper, the authors studied the variation of dimension under holomorphic motions and proved the following results. We denote by $\dimH, \dimP, \dimMl, \dimMu$ the Hausdorff, packing, lower Minkowski and upper Minkowski dimensions respectively. In the first result $\dim$ denotes either $\dimMu$ or $\dimP$.

\begin{theorem}[Fuhrer--Ransford--Younsi]
\label{T:FRY1}
Let $\dim$ denote either of $\dimP$ or $\dimMu$, and let $f: \mD \times A \to \mC$ be a holomorphic motion of a set $A \subset \mC$, which we assume to be bounded in the case $\dim=\dimMu$. Then either $\dim(A_\lambda)=0$ for all $\lambda \in \mD$, or $\lambda \mapsto 1/\operatorname{dim}(A_\lambda)$ is an inf-harmonic function on $\mD$. In particular, $\lambda\mapsto\operatorname{dim}(A_\lambda)$ is a subharmonic function on $\mD$.
\end{theorem}

The authors also proved that, if $\dimH(A_0)>0$, then  the function $\lambda \mapsto 1/\dimH(A_\lambda)$ is the supremum of a family of inf-harmonic functions, and asked whether it is in fact inf-harmonic; see \cite[Question 11.1]{FRY}. They  remarked that 
a positive answer would imply that 
$\lambda \mapsto \dimH(A_\lambda)$ is a subharmonic function, 
and that even this weaker statement was still an open problem.

The goal of this note is to answer both of these questions in the negative, not only for $\dimH$, but also for $\dimMl$. Our main result is as follows.

\begin{theorem}
\label{T:main}
There exist a compact set $A\subset \mC$ and a holomorphic motion $f:\mD\times A\to\mC$ such that, writing
$M(\lambda):=\max\{|\operatorname{Re}\lambda|,|\operatorname{Im}\lambda|\}$,
we have
\begin{equation}\label{eq:hausdorff-main}
 \dimH(A_\lambda)=\dimMl(A_\lambda)=\frac{1}{3+M(\lambda)}
 \qquad(\lambda\in\mD),
\end{equation}
and
\begin{equation}\label{eq:packing-main}
 \dimP(A_\lambda)=\dimMu(A_\lambda)=\frac{1}{3-M(\lambda)}
 \qquad(\lambda\in\mD).
\end{equation}
In particular, as both $\dimH(A_\lambda)$ and $\dimMl(A_\lambda)$
attain a strict maximum at $\lambda=0$,  neither is  subharmonic,
and the reciprocal of neither  is inf-harmonic.
\end{theorem}

We give the construction and its verification in the following sections,
and we make some general comments on the example at the end.

\section{Construction}\label{S:constr}
Set
\[
 N_0:=0,\qquad N_k:=k!\quad(k\ge1),
\]
and define
\[
 \sigma_n:=i^{\,k-1}\qquad\text{if }N_{k-1}<n\le N_k,
\]
so the sequence $(\sigma_n)$ contains the values $1,i,-1,-i$ alternating cyclically in constant blocks. We also define 
\[
 S_0:=0,\qquad S_n:=\sum_{j=1}^n\sigma_j\quad(n\ge1).
\]

\begin{lemma}
\label{L:limit}
For each $\lambda \in \mC$, we have
\begin{equation}\label{eq:limsup}
 \limsup_{n\to\infty}\operatorname{Re}\left(\lambda\frac{S_n}{n}\right)
 =\max\{|\operatorname{Re}\lambda|,|\operatorname{Im}\lambda|\}
\end{equation}
and
\begin{equation}\label{eq:liminf}
 \liminf_{n\to\infty}\operatorname{Re}\left(\lambda\frac{S_n}{n}\right)
 =-\max\{|\operatorname{Re}\lambda|,|\operatorname{Im}\lambda|\}.
\end{equation}

\end{lemma}

\begin{proof}
It suffices to prove (\ref{eq:limsup}), since (\ref{eq:liminf}) then follows by replacing $\lambda$ with $-\lambda$. 

Note that for $k \geq 2$ we have

\begin{equation}\label{eq:block-end}
\begin{aligned}
 \left|\frac{S_{N_k}}{N_k}-i^{\,k-1}\right|
 &=\left|\frac{S_{N_{k-1}}-N_{k-1}i^{\,k-1}}{N_k}\right|\\
 &\le\frac{2N_{k-1}}{N_k}=\frac2k.
\end{aligned}
\end{equation}
In particular each of $1,i,-1,-i$ is the  limit of some subsequence of $S_n/n$, from which it easily follows that
\[
\max\{|\operatorname{Re}\lambda|,|\operatorname{Im}\lambda|\} \leq \limsup_{n\to\infty}\operatorname{Re}\left(\lambda\frac{S_n}{n}\right).
\]

For the reverse inequality, note that every $\sigma_j$ belongs to $\{1,i,-1,-i\}$, and thus
\[
\operatorname{Re}\left(\lambda\frac{S_n}{n}\right) = \frac{1}{n} \sum_{j=1}^n \operatorname{Re}\left(\lambda \sigma_j\right) \leq \max\{|\operatorname{Re}\lambda|,|\operatorname{Im}\lambda|\},
\]
as required.
\end{proof}

Now, for $n \geq 0$ and $\lambda \in \mD$, define
\[
 R_n(\lambda):=2^{-3n-\lambda S_n}.
\]

Let $\Omega:=\{0,1\}^\mN$, equipped with the product topology. For $\lambda \in \mD$ and $\varepsilon=(\varepsilon_n)_{n\ge1}\in\Omega$, define
\[
 \pi_\lambda(\varepsilon):=\sum_{n=1}^{\infty}\varepsilon_nR_n(\lambda).
\]

\begin{lemma}
\label{L:motion}
\begin{enumerate}[\normalfont(i)]
\item For each $\lambda \in \mD$, the function $\pi_\lambda:\Omega \to \mC$ is continuous and injective.
\item For each $\varepsilon \in \Omega$, the function $\lambda \mapsto \pi_\lambda(\varepsilon)$ is holomorphic on $\mD$.
\end{enumerate}
\end{lemma}

\begin{proof}
Note that $R_0=1$. For $n \in \mN$ and $\lambda \in \mD$, we have
\[
 \frac{R_n(\lambda)}{R_{n-1}(\lambda)}
 =2^{-3-\lambda\sigma_n}.
\]
Since $|\sigma_n|=1$ and $|\lambda|<1$, this implies that
\begin{equation}\label{eq:ratio}
 \frac{1}{16}
 \le \left|\frac{R_n(\lambda)}{R_{n-1}(\lambda)}\right|
 \le \frac{1}{4}.
\end{equation}
In particular $|R_n(\lambda)|\le 4^{-n}$. It follows that the series
\[
 \pi_\lambda(\varepsilon)
 :=\sum_{n=1}^{\infty}\varepsilon_nR_n(\lambda)
\]
converges uniformly in both $\lambda \in \mD$ and
$\varepsilon \in \Omega$. Thus the function
$\pi_\lambda:\Omega \to \mC$ is continuous, and, for each fixed
$\varepsilon \in \Omega$, the function
$\lambda \mapsto \pi_\lambda(\varepsilon)$ is holomorphic on $\mD$.

It remains to show that $\pi_\lambda$ is injective.
Let $\varepsilon,\eta \in \Omega$ with $\varepsilon \neq \eta$,
and suppose that they first differ at position $m$.
By \eqref{eq:ratio}, we have
\begin{equation}
\label{eq1}
 \sum_{n>m}|R_n(\lambda)|
 \le |R_m(\lambda)|\sum_{j=1}^{\infty}4^{-j}
 =\frac{1}{3}|R_m(\lambda)|.
 \end{equation}
As $|\varepsilon_m-\eta_m|=1$ and
$|\varepsilon_n-\eta_n|\le 1$ for every $n>m$, we obtain
\[
\begin{aligned}
 |\pi_\lambda(\varepsilon)-\pi_\lambda(\eta)|
 &= \left|
     (\varepsilon_m-\eta_m)R_m(\lambda)
     +\sum_{n>m}(\varepsilon_n-\eta_n)R_n(\lambda)
    \right| \\
 &\ge |R_m(\lambda)|-\sum_{n>m}|R_n(\lambda)| \ge \frac{2}{3}|R_m(\lambda)|. 
\end{aligned}
\]
Therefore $\pi_\lambda(\varepsilon) \neq \pi_\lambda(\eta)$ and $\pi_\lambda:\Omega \to \mC$ is injective, as required.
\end{proof}

Now define
\[
 A:=\pi_0(\Omega)
 =\left\{\sum_{n=1}^{\infty}\varepsilon_n8^{-n}:
       \varepsilon_n\in\{0,1\}\right\}.
\]
Geometrically, $A$ is a Cantor set in $\mR$ constructed as follows. Starting with $[0,1/7]$, we retain the leftmost and rightmost one-eighths of the interval, removing the open middle three-quarters. We repeat this procedure in each remaining interval. The set $A$ is the intersection of the resulting decreasing sequence of compact sets.

Define $f:\mD \times A \to \mC$ by
\[
f\bigl(\lambda,\pi_0(\varepsilon)\bigr):=\pi_\lambda(\varepsilon).
\]

By Lemma \ref{L:motion}, the map $f$ defines a holomorphic motion of $A$. Moreover, for each $\lambda\in\mD$, the set $A_\lambda=\pi_\lambda(\Omega)$ is compact, since $\Omega$ is compact and $\pi_\lambda$ is continuous.

\section{Hausdorff dimension}\label{S:hausdorff}

In this section we prove the desired formula for the Hausdorff dimension of the sets $A_\lambda$ from Theorem \ref{T:main}, namely

\begin{equation}\label{eq:hausdorff}
 \dimH(A_\lambda)=\frac{1}{3+M(\lambda)} \qquad(\lambda\in\mD),
\end{equation}
where $M(\lambda):=\max\{|\operatorname{Re}\lambda|,|\operatorname{Im}\lambda|\}$.

\begin{lemma}
\label{L:hausdorff}
For each $\lambda \in \mD$, we have
\[
\dimH(A_\lambda)=\liminf_{n\to\infty}\frac{n\log2}{-\log|R_n(\lambda)|}.
\]
\end{lemma}

Equation (\ref{eq:hausdorff}) easily follows from Lemma \ref{L:hausdorff}. Indeed, by definition of $R_n(\lambda)$, we have
\[
 -\log|R_n(\lambda)|=(3n+\operatorname{Re}(\lambda S_n))\log2.
\]
Combining this with \eqref{eq:limsup}, we obtain
\[
\begin{aligned}
 \dimH(A_\lambda)
 &=\liminf_{n\to\infty}\frac1{3+\operatorname{Re}(\lambda S_n/n)}\\
 &=\frac1{3+\limsup_{n\to\infty}\operatorname{Re}(\lambda S_n/n)}\\
 &=\frac1{3+\max\{|\operatorname{Re}\lambda|,|\operatorname{Im}\lambda|\}}\\
 &=\frac1{3+M(\lambda)},
\end{aligned}
\]
as required.

For the proof of Lemma \ref{L:hausdorff}, we introduce the following definition.

\begin{definition}
Fix $\lambda \in \mD$. A \textit{cylinder} of length $n$ is the subset of $A_\lambda$ obtained by fixing the first $n$ binary digits 
$\varepsilon_1,\dots,\varepsilon_n$ in 
$$
 \pi_\lambda(\varepsilon):=\sum_{j=1}^{\infty}\varepsilon_jR_j(\lambda).
$$
\end{definition}
We will need the following properties of cylinders.
\begin{proposition}
\label{P:cylinders}
Fix $\lambda \in \mD$ and write $\ell_n:=|R_n(\lambda)|$.
\begin{enumerate}[\normalfont(i)]
\item Each cylinder of length $n$ has diameter at most $\ell_n/3$.
\item The distance between two distinct cylinders of length $n$ is at least $2\ell_n/3$.
\end{enumerate}
\end{proposition}

\begin{proof}
By definition, the diameter of any cylinder of length $n$ is at most
\[
\sum_{m>n} \ell_m \leq \frac{1}{3}\ell_n,
\]
by \eqref{eq1}. This proves (i).

For (ii), let  $\pi_\lambda(\varepsilon)$ and $\pi_\lambda(\eta)$ 
belong to distinct cylinders of length $n$. Then $\varepsilon$ and $\eta$ first differ at some position $m \leq n$. As in the proof of Lemma \ref{L:motion}, we obtain
\[
|\pi_\lambda(\varepsilon)-\pi_\lambda(\eta)|
\ge |R_m(\lambda)|-\sum_{j>m}|R_j(\lambda)|
\ge \frac23\ell_m
\ge \frac23\ell_n,
\]
as required.
\end{proof}

\begin{proof}[Proof of Lemma \ref{L:hausdorff}]
Fix $\lambda\in\mD$, and again write $\ell_n=|R_n(\lambda)|$. Let 
\[
D:=\liminf_{n\to\infty}\frac{n\log2}{-\log\ell_n}.
\]
We have to show that $\dimH(A_\lambda)=D$.

First note that there are $2^n$ cylinders of length $n$, 
and that each has diameter at most $\ell_n/3$ by Proposition~\ref{P:cylinders}. Fix $s>D$, and choose $t$ with $D<t<s$. There are infinitely many $n$ for which
\[
 \frac{n\log2}{-\log\ell_n}<t.
\]
Along this subsequence,
\[
 2^n\left(\frac{\ell_n}{3}\right)^s
 \le3^{-s}\ell_n^{\,s-t}\longrightarrow0.
\]
The cylinder covers therefore give 
the estimate $\dimH(A_\lambda)\le s$.
As this holds for each $s>D$, we conclude that
$\dimH(A_\lambda)\le D$.

For the reverse inequality, we consider the push-forward under $\pi_\lambda$ of the product probability measure on $\Omega$ that gives probability $1/2$ to each digit. The resulting probability measure $\mu$ on $A_\lambda$ assigns mass $2^{-n}$ to every cylinder of length $n$.

Fix $0<s<D$. For all sufficiently large $n$, we have
\begin{equation}\label{eq:mass-scale}
 2^{-n}\le\ell_n^s.
\end{equation}
For a sufficiently small radius $r>0$, choose $n$ such that
\begin{equation}\label{eq:radius-scale}
 \frac{\ell_n}{3}\le r<\frac{\ell_{n-1}}{3}.
\end{equation}
By Proposition \ref{P:cylinders}, the distance between two different cylinders of length $n-1$ is at least
\[
 \frac23\ell_{n-1}>2r.
\]
A ball of radius $r$ thus meets at most one cylinder of length $n-1$.
Combining \eqref{eq:mass-scale} and \eqref{eq:radius-scale}, we obtain
\begin{equation}\label{eq:mass-bound}
 \mu(B(z,r))\le2^{-(n-1)}=2\,2^{-n}\le2\ell_n^s\le2(3r)^s.
\end{equation}
This holds for all centers $z\in\mC$ and all sufficiently small $r$. It follows from the Mass Distribution Principle \cite[Lemma~1.2.8]{BP} that $\dimH(A_\lambda) \geq s$. Since this holds for every $0<s<D$, we deduce that $\dimH(A_\lambda) \geq D$, as required.
\end{proof}

\section{Lower and upper Minkowski dimensions}\label{S:minkowski}

In this section we prove the  formulas for the lower and upper Minkowski dimensions of the sets $A_\lambda$ from Theorem \ref{T:main}, namely

\begin{equation}\label{eq:lowerM}
 \dimMl(A_\lambda)=\frac{1}{3+M(\lambda)}
 \qquad(\lambda\in\mD)
\end{equation}
and
\begin{equation}\label{eq:upperM}
 \dimMu(A_\lambda)=\frac{1}{3-M(\lambda)}
 \qquad(\lambda\in\mD),
\end{equation}
where $M(\lambda):=\max\{|\operatorname{Re}\lambda|,|\operatorname{Im}\lambda|\}$.

\begin{lemma}
\label{L:minkowski}
For each $\lambda \in \mD$, we have
\[
 \dimMl(A_\lambda)=\liminf_{n\to\infty}\frac{n\log2}{-\log\ell_n}
\]
and
\[
 \dimMu(A_\lambda)=\limsup_{n\to\infty}\frac{n\log2}{-\log\ell_n},
\]
where $\ell_n=|R_n(\lambda)|$.
\end{lemma}

Equations \eqref{eq:lowerM} and  \eqref{eq:upperM} easily follow from Lemma \ref{L:minkowski}, using the relation
\[
 \frac{n\log2}{-\log\ell_n}
 =\frac1{3+\operatorname{Re}(\lambda S_n/n)}
\]
and Lemma \ref{L:limit}.

\begin{proof}[Proof of Lemma \ref{L:minkowski}]
Fix $\lambda\in\mD$, and again set $\ell_n:=|R_n(\lambda)|$. For a bounded set $E\subset\mC$ and $\delta>0$, let $N_\delta(E)$ be the smallest number of sets of diameter at most $\delta$ needed to cover $E$. Recall that the lower and upper Minkowski dimensions are defined by
\[
 \dimMl(E):=\liminf_{\delta\downarrow0}
       \frac{\log N_\delta(E)}{\log(1/\delta)},\qquad
 \dimMu(E):=\limsup_{\delta\downarrow0}
       \frac{\log N_\delta(E)}{\log(1/\delta)}.
\]

For any sufficiently small $\delta>0$, choose $n$ such that
\begin{equation}\label{eq:delta-scale}
 \frac{\ell_n}{3}\le\delta<\frac{\ell_{n-1}}3.
\end{equation}
The cylinders of length $n$ provide a cover of $A_\lambda$, whereas a set of diameter at most $\delta$ cannot meet two cylinders of length $n-1$, by Proposition \ref{P:cylinders}. Consequently,
\begin{equation}\label{eq:covering-bounds}
 2^{n-1}\le N_\delta(A_\lambda)\le2^n.
\end{equation}
The bounds $1/16\le\ell_n/\ell_{n-1}\le1/4$ from \eqref{eq:ratio} imply that, for the choice
of $n$ in \eqref{eq:delta-scale},
\[
 -\log\delta=-\log\ell_n+O(1),\qquad
 \log N_\delta(A_\lambda)=n\log2+O(1).
\]
The constants are independent of $n$ and $\delta$. It follows that
\[
 \dimMl(A_\lambda)=\liminf_{n\to\infty}\frac{n\log2}{-\log\ell_n}
\]
and
\[
 \dimMu(A_\lambda)=\limsup_{n\to\infty}\frac{n\log2}{-\log\ell_n},
\]
as required.
\end{proof}

\section{Packing dimension}\label{S:packing}

In this section we prove the desired formula for the packing dimension of the sets $A_\lambda$ from Theorem \ref{T:main}, namely

\begin{equation}\label{eq:packing}
 \dimP(A_\lambda)=\frac{1}{3-M(\lambda)} \qquad(\lambda\in\mD),
\end{equation}
where $M(\lambda):=\max\{|\operatorname{Re}\lambda|,|\operatorname{Im}\lambda|\}$.

\begin{lemma}
\label{L:packing}
We have
\[
\dimP(A_\lambda)=\dimMu(A_\lambda)  \qquad(\lambda\in\mD).
\]
\end{lemma}

Equation \eqref{eq:packing} follows directly from Lemma \ref{L:packing} and \eqref{eq:upperM}.

\begin{proof}[Proof of Lemma \ref{L:packing}]
We use the following characterization of packing dimension due to Tricot:
\begin{equation}\label{eq:packing-cover}
 \dimP(E)=\inf\left\{\sup_{j\ge1}\dimMu(E_j):
                E=\bigcup_{j\ge1}E_j\right\},
\end{equation}
where the infimum is over countable covers by bounded subsets $E_j$ of $E$; see \cite[Proposition 2.3]{FRY}.

Fix $\lambda\in\mD$, write $\ell_n:=|R_n(\lambda)|$, and put
\[
 p:=\dimMu(A_\lambda)=\frac1{3-M(\lambda)}.
\]
Taking the cover consisting of $A_\lambda$ alone in \eqref{eq:packing-cover} gives $\dimP(A_\lambda)\le p$.

For the reverse inequality, let $\mu$ be the probability measure on $A_\lambda$ which assigns mass $2^{-n}$ to every cylinder of length $n$. Write $\mu^*$ for the associated outer measure.

Suppose that $E\subset A_\lambda$ satisfies $\mu^*(E)>0$. 
Set $a:=\mu^*(E)$. Let $m_n$ be the number of cylinders of length $n$ intersecting $E$. The union of all such cylinders contains $E$, and hence
\[
 a\le m_n2^{-n}.
\]
By Proposition \ref{P:cylinders}, distinct cylinders of length $n$ are separated by at least $2\ell_n/3$. In particular a set of diameter at most $\ell_n/3$ cannot meet two of them. It follows that
\[
 N_{\ell_n/3}(E)\ge m_n\ge a\,2^n.
\]
Taking a limsup along these scales and using Lemma \ref{L:minkowski} gives
\[
 \dimMu(E)\ge\limsup_{n\to\infty}
       \frac{n\log2+\log a}{-\log\ell_n+\log3}=p.
\]
Monotonicity gives the reverse inequality, because $E\subset A_\lambda$. We have therefore shown that,
whenever $E\subset A_\lambda$,
\begin{equation}\label{eq:positive-mass}
 \mu^*(E)>0\quad\Longrightarrow\quad\dimMu(E)=p.
\end{equation}
Now, for any countable cover $A_\lambda=\bigcup_{j\ge1}E_j$, at least one $E_j$ has positive $\mu$-outer measure. By \eqref{eq:positive-mass}, that set
has upper Minkowski dimension $p$. Thus every cover in \eqref{eq:packing-cover} satisfies
\[
 \sup_{j\ge1}\dimMu(E_j)\ge p,
\]
from which it follows that $\dimP(A_\lambda)\ge p$, as required.
\end{proof}

\section{Concluding remarks}

Summarizing the results of Sections \ref{S:hausdorff}, \ref{S:minkowski} and  \ref{S:packing}, we have established Equations \eqref{eq:hausdorff-main} and \eqref{eq:packing-main} in Theorem \ref{T:main}. The last sentence in the statement of the theorem is straightforward: the Hausdorff and lower Minkowski dimension functions both attain a strict maximum at $0$, so, by the maximum principle, neither function 
can be subharmonic, and hence neither reciprocal is inf-harmonic.
Thus Theorem~\ref{T:main} is proved.

We conclude by mentioning that, for  the holomorphic motion $f:\mD\times A\to\mC$ of Theorem \ref{T:main}, the reciprocals of packing dimension and upper Minkowski dimension are inf-harmonic, which is consistent with Theorem \ref{T:FRY1}. Indeed, for $\lambda:=x+iy \in \mD$, we have
\[
\frac1{\dimP(A_\lambda)}=\frac1{\dimMu(A_\lambda)}=3-M(\lambda)
 =\min\{3+x,3-x,3+y,3-y\},
\]
and each function in this minimum is positive and harmonic on $\mD$. 

The example is also consistent with Theorem~1.6 of~\cite{FRY}, which asserts
that $1/\dimH(A_\lambda)-1/2$ is a supremum of inf-harmonic functions.
Indeed, again writing $\lambda=x+iy$, we have
\[
\frac1{\dimH(A_\lambda)}-\frac12
=\max\left\{
\frac52+x,\frac52-x,\frac52+y,\frac52-y
\right\},
\]
which is a maximum of four positive harmonic functions on $\DD$.

\section*{AI tools declaration}

In preparing this work, the authors used ChatGPT to develop and explore strategies for the proof of Theorem~\ref{T:main}. The final mathematical  text was written entirely by the authors, who thoroughly corrected, verified, and streamlined the proposed strategy. The authors accept full responsibility for the accuracy and integrity of the final results presented in this paper.



\begin{thebibliography}{9}
\bibitem{BP}
C.~J. Bishop and Y.~Peres,
\emph{Fractals in Probability and Analysis},
Cambridge Studies in Advanced Mathematics, vol.~162,
Cambridge University Press, Cambridge, 2017.

\bibitem{FRY}
A.~Fuhrer, T.~Ransford and M.~Younsi,
Holomorphic motions, dimension, area and quasiconformal mappings,
\emph{J. Math. Pures Appl.} 177 (2023), 455--483.



\end{thebibliography}
\end{document}